\documentclass{amsart}
\usepackage{amsmath} 
\usepackage{amssymb}
\usepackage{mathrsfs}

\newtheorem{theorem}{Theorem}[section] 
\newtheorem{claim}[theorem]{Claim}
\newtheorem{lemma}[theorem]{Lemma} 
 
\newtheorem{observation}[theorem]{Observation} 
 
\newtheorem{conclusion}[theorem]{Conclusion}

\theoremstyle{definition}
\newtheorem{definition}[theorem]{Definition}

\newtheorem{defclaim}[theorem]{Definition/Claim}

\theoremstyle{remark}
\newtheorem{remark}[theorem]{Remark}
\newtheorem{notation}[theorem]{Notation}
\newtheorem{question}[theorem]{Question}

\newcommand{\rest}{{\restriction}}
 
\newcommand{\PC}{{\rm PC}} 
 
\newcommand{\tr}{{\rm tr}} 
\newcommand{\tp}{{\rm tp}} 
\newcommand{\Av}{{\rm Av}} 
\newcommand{\EC}{{\rm EC}} 
\newcommand{\OB}{{\rm OB}} 
\newcommand{\EM}{{\rm EM}} 
\newcommand{\Car}{{\rm Car}} 
\newcommand{\Min}{{\rm Min}} 

\newcommand{\stat}{{\rm stat}}

\newcommand{\then}{{\underline{then}}}
\newcommand{\when}{{\underline{when}}}
\newcommand{\Then}{{\underline{Then}}}

\newcommand{\If}{{\underline{if}}}
\newcommand{\Iff}{{\underline{iff}}}
\newcommand{\mn}{{\medskip\noindent}}
\newcommand{\sn}{{\smallskip\noindent}}

\newcommand{\bbL}{{\mathbb L}}

\newcommand{\gk}{{\mathfrak k}}

\newcommand{\bbP}{{\mathbb P}}
\newcommand{\cP}{{\mathscr P}}

\newcommand{\cU}{{\mathscr U}}

\newcommand{\cf}{{\rm cf}}

\newcount\skewfactor
\def\mathunderaccent#1#2 {\let\theaccent#1\skewfactor#2
\mathpalette\putaccentunder}
\def\putaccentunder#1#2{\oalign{$#1#2$\crcr\hidewidth
\vbox to.2ex{\hbox{$#1\skew\skewfactor\theaccent{}$}\vss}\hidewidth}}

\newenvironment{PROOF}[2][\proofname.]
   {\begin{proof}[#1]\renewcommand{\qedsymbol}{\textsquare$_{\rm #2}$}}
   {\end{proof}}

\begin{document}

\title {When first order $T$ has limit models}
\author {Saharon Shelah}
\address{Einstein Institute of Mathematics\\
Edmond J. Safra Campus, Givat Ram\\
The Hebrew University of Jerusalem\\
Jerusalem, 91904, Israel\\
 and \\
 Department of Mathematics\\
 Hill Center - Busch Campus \\ 
 Rutgers, The State University of New Jersey \\
 110 Frelinghuysen Road \\
 Piscataway, NJ 08854-8019 USA}
\email{shelah@math.huji.ac.il}
\urladdr{http://shelah.logic.at}
\thanks {The author would like to thank the Israel Science Foundation
 for partial support of this research (Grant No. 710/07).
I would like to thank Alice Leonhardt 
for the beautiful typing.  First Typed - 04/June/23.  Paper 868 in the
 author list of publications.}



\date{March 8, 2013}



\begin{abstract}
We to a large extent sort out when does a (first order
complete theory) $T$ have a superlimit model in a cardinal $\lambda$.  
Also we deal with related notions of being limit.
\end{abstract}

\maketitle
\numberwithin{equation}{section}
\setcounter{section}{-1}
\newpage

\section*{Anotated Content}
\bigskip

\noindent
\S0 \quad Introduction, pg.\pageref{Introduction}
\mn
\begin{enumerate}
\item[${{}}$]    [We give background and the basic definitions.  We
then present existence results for stable $T$ which have models which
are saturated or closed to being saturated.]
\end{enumerate}
\bigskip

\noindent
\S1 \quad On countable superstable not $\aleph_0$-stable,
pg.\pageref{Onsuperstable} 
\mn
\begin{enumerate}
\item[${{}}$]    [Consistently $2^{\aleph_1} \ge \aleph_2$ and some such
(complete first order) $T$ has a superlimit (non-saturated) model
of cardinality $\aleph_1$.  This shows that we cannot prove a
non-existence result fully complementary to Lemma \ref{y.7}.]
\end{enumerate}
\bigskip

\noindent
\S2 \quad A strictly stable consistent example, pg.\pageref{Astrictly}
\mn
\begin{enumerate}
\item[${{}}$]   [Consistently $\aleph_1 < 2^{\aleph_0}$ and 
some countable stable not superstable $T$, has a (non-saturated) model
of cardinality $\aleph_1$ which satisfies some relatives of being
superlimit.]
\end{enumerate}
\bigskip

\noindent
\S3 \quad On the non-existence of limit models, pg.\pageref{On non}
\mn
\begin{enumerate}
\item[${{}}$]   [The proofs here are in ZFC.  If $T$ is unstable it
has no superlimit models of cardinality $\lambda$ when $\lambda 
\ge \aleph_1 + |T|$.  For
unsuperstable $T$ we have similar results but with ``few" exceptional
cardinals $\lambda$ on which we do not know: $\lambda < \lambda^{\aleph_0}$
which are $< \beth_\omega$.  Lastly, if $T$ is superstable and
$\lambda \ge |T| + 2^{|T|}$ then $T$ has a superlimit model of
cardinality $\lambda$ iff $|D(T)| \le \lambda$ iff $T$ has a
saturated model.  Lastly, we get weaker results on weaker relatives of
superlimit.] 
\end{enumerate}
\newpage

\section {Introduction} \label{Introduction}
\bigskip

\noindent
\S (0A) \quad Background and Content

Recall that (\cite[Ch.III]{Sh:c}).
If $T$ is (first order complete and) superstable then for $\lambda \ge
2^{|T|},T$ has a saturated model $M$ of cardinality $\lambda$ and
moreover
\mn
\begin{enumerate}
\item[$(*)$]  if $\langle M_\alpha:\alpha < \delta\rangle$ is
  $\prec$-increasing, $\delta$ a limit ordinal $< \lambda^+$ and
  $\alpha < \delta \Rightarrow M_\alpha \cong M$ then $\cup
  \{M_\alpha:\alpha < \delta\}$ is isomorphic to $M$.
\end{enumerate}
\mn
When investigating categoricity of an a.e.c. (abstract elementary
classes) $\gk = (K_{\gk},\le_{\gk})$, the following property turns out
to be central: $M$ is $\le_{\gk}$-universal model of cardinality
$\lambda$ with the property $(*)$ above (called superlimit) - 
possibly with addition parameter $\kappa = \cf(\kappa) \le
\lambda$ (or stationary $S \subseteq \lambda^+$); we also consider
some relatives, mainly limit, weakly limit and strongly limit.  
Those notions were suggested for a.e.c. in \cite[3.1]{Sh:88} or see 
the revised version \cite[3.3]{Sh:88r} and see \cite{Sh:h}
\underline{or} here in \ref{y.5}.  But though coming from
investigating non-elementary classes, they are meaningful for
elementary classes and here we try to investigate them for elementary
classes.  

Recall that for a first order complete $T$, 
we know $\{\lambda:T$ has a saturated
model of $T$ of cardinality $\lambda\}$, that is, it is $\{\lambda:\lambda^{<
\lambda} \ge |D(T)|$ or $T$ is stable in $\lambda\}$, on the
definitions of $D(T)$ and
other notions see \S(0B) below.  What if we replace saturated by 
superlimit (or some relative)?   Let EC$_\lambda(T)$ be the class
of models $M$ of $T$ of cardinality $\lambda$.

If there is a saturated $M \in \EC_\lambda(T)$ we have
considerable knoweldge on the existence of limit model for cardinal
$\lambda$, this was as mentioned in \cite[3.6]{Sh:88r} by \cite{Sh:c},
see \ref{y.7}(1),(2).  E.g. for superstable $T$ in $\lambda \ge 
2^{|T|}$ there is a superlimit model 
(the saturated one).  It seems a natural question on 
\cite[3.6]{Sh:88r} whether it exhausts the
possibilities of $(\lambda,*)$-superlimit and
$(\lambda,\kappa)$-superlimit models for elementary classes.

Clearly the cases of the existence of
such models of a (first order complete) theory $T$ where there are no
saturated (or special) models are rare, because
even the weakest version of Definition \cite[3.1]{Sh:88} =
\cite[3.3]{Sh:88r} or here Definitino \ref{y.5}
for $\lambda$ implies that $T$ has a universal model
of cardinality $\lambda$, which is rare (see Kojman Shelah
\cite{KjSh:409} which includes earlier history and recently Djamonza
\cite{Dj04}).

So the main question seems to be whether there are such cases at all.  
We naturally look at some of the previous cases of
consistency of the existence of a universal model (for $\lambda < \lambda^{<
\lambda}$), i.e., those for $\lambda = \aleph_1$.

E.g. a sufficient condition for some versions is the existence of $T'
\supseteq T$ of cardinality $\lambda$ such that $\PC(T',T)$ is
categorical in $\lambda$, see \ref{y.2}(3).  By \cite{Sh:100} we have
consistency results for such $T_1$ so naturally we 
first deal with the consistency results from \cite{Sh:100}.
In \S1 we deal with the case of the countable superstable $T_0$ from
\cite{Sh:100} which is
not $\aleph_0$-stable. By \cite{Sh:100} consistently $\aleph_1 <
2^{\aleph_0}$ and for some $T'_0 \supseteq T_0$ of cardinality
$\aleph_1$, PC$(T'_0,T_0)$ is categorical in $\aleph_1$.  We use this to
get the consistency of ``$T_0$ has a superlimit model of cardinality
$\aleph_1$ and $\aleph_1 < 2^{\aleph_0}$".

In \S2 we prove that 
for some stable not superstable countable $T_1$ we have a parallel
but weaker result.  
We relook at the old consistency results of ``some
$\PC(T'_1,T_1),|T'_1| = \aleph_1 > |T_1|$, is 
categorical in $\aleph_1$" from \cite{Sh:100}.  From this we
deduce that in this universe, $T_1$ has a strongly
$(\aleph_1,\aleph_0)$-limit model.

It is a reasonable thought that we can similarly have a consistency
result on the theory of linear order, but this is still not clear.

In \S3 we show that if $T$ has a superlimit model in $\lambda \ge |T|
+ \aleph_1$ then $T$ is stable and $T$ is superstable except
possibly under some severe restrictions on the cardinal $\lambda$
(i.e., $\lambda < \beth_\omega$ and $\lambda < \lambda^{\aleph_0}$).
We then prove some restrictions on the existence of some (weaker)
relatives.

Summing up our results on the strongest notion, superlimit, by
\ref{nlm.0.0} + \ref{nlm.0.1} we have:
\begin{conclusion}
\label{n0.1}  
Assume $\lambda \ge |T| +
\beth_\omega$.  \Then \, $T$ has a superlimit model of cardinality
$\lambda$ iff $T$ is superstable and $\lambda \ge |D(T)|$.
\end{conclusion}

\noindent
In subsequent work we shall show that for some unstable $T$
(e.g. the theory of linear orders), if 
$\lambda = \lambda^{< \lambda} > \kappa = \cf(\kappa)$, \then \, 
$T$ has a medium $(\lambda,\kappa)$-limit model,
whereas if $T$ has the independence property even weak
$(\lambda,\kappa)$-limit models do not exist; see \cite{Sh:877} and
more in \cite{Sh:900}, \cite{Sh:906}, \cite{Sh:950}, \cite{Sh:F1054}.

We thank Alex Usvyatsov for urging us to resolve the question of the
superlimit case and John Baldwin for comments and complaints.
\bigskip

\noindent
\S (0B) \quad Basic Definitions

\begin{notation}
\label{y.1}
1) Let $T$ denote a complete first order theory which has infinite
   models but $T_1,T'$, etc. are not necessarily complete.

\noindent
2) Let $M,N$ denote models, $|M|$ the universe of $M$ and $\|M\|$ its
   cardinality and $M \prec N$ means $M$ is an elementary submodel of
   $N$.

\noindent
3) Let $\tau_T = \tau(T),\tau_M = \tau(M)$ be the 
vocabulary of $T,M$ respectively.

\noindent
4) Let $M \models ``\varphi[\bar a]^{\If(\stat)}"$ means that the
   model $M$ satisfies $\varphi[\bar a]$ iff the statement $\stat$ 
is true (or is 1 rather than 0)).
\end{notation}

\begin{definition}
\label{y1d}
1) For $\bar a \in {}^{\omega >}|M|$ and $B \subseteq M$ let $\tp(\bar
   a,B,M) = \{\varphi(\bar x,\bar b):\varphi = \varphi(\bar x,\bar y)
   \in \bbL(\tau_M),\bar b \in {}^{\ell g(\bar y )}B$ and $M \models
   \varphi[\bar a,\bar b]\}$.

\noindent
2) Let $D(T) = \{\tp(\bar a,\emptyset,M):M$ a model of $T$ and $\bar
   a$ a finite sequence from $M\}$.

\noindent
3) If $A \subseteq M$ then $\bold S^m(A,M) = \{\tp(\bar a,A,N):M \prec
   N$ and $\bar a \in {}^mN\}$, if $m=1$ we may omit it.

\noindent
4) A model $M$ is $\lambda$-saturated when: if $A \subseteq M,|A| <
   \lambda$ and $p \in \bold S(A,M)$ then $p$ is realized by some $a
   \in M$, i.e. $p \subseteq \tp(a,A,M)$; if $\lambda = \|M\|$ we may
   omit it.

\noindent
5) A model $M$ is special \when \, letting $\lambda = \|M\|$, there is
   an increasing sequence $\langle \lambda_i:i < \cf(\lambda)\rangle$
   of cardinals with limit $\lambda$ and a $\prec$-increasing sequence
   $\langle M_i:i < \cf(\lambda)\rangle$ of models with union $M$ such that
   $M_{i+1}$ is $\lambda_i$-saturated of cardinality $\lambda_{i+1}$
   for $i < \cf(\lambda)$.
\end{definition}

\begin{definition}
\label{y.2}  
1) For any $T$ let EC$(T) = \{M:M$ is a $\tau_T$-model of $T\}$.

\noindent
2) EC$_\lambda(T) = \{M \in \text{ EC}(T):M$ is of cardinality
   $\lambda\}$.

\noindent
3) For $T \subseteq T'$ let

\[
\text{PC}(T',T) = \{M \restriction \tau_T:M \text{ is model of } T'\}
\]

\[
\text{PC}_\lambda(T',T) = \{M \in \text{ PC}(T',T):M \text{ is of
cardinality } \lambda\}.
\]
\mn
4) We say $M$ is $\lambda$-universal for $T_1$ when it is a model of
$T_1$ and every $N \in \text{ EC}_\lambda(T)$ can be elementarily 
embedded into $M$; if $T_1 = \text{Th}(M)$ we may omit it.

\noindent
5) We say $M \in \text{ EC}(T)$ is universal \when \, it is
$\lambda$-universal for $\lambda = \|M\|$.
\end{definition}

\noindent
We are here mainly interested in
\begin{definition}
\label{y.3}  
Given $T$ and $M \in \text{ EC}_\lambda(T)$ we say that 
$M$ is a superlimit or $\lambda$-superlimit model \when \,: $M$ is
universal and if $\delta < \lambda^+$ is a limit
ordinal, $\langle M_\alpha:\alpha \le \delta \rangle$ is
$\prec$-increasing continuous, and $M_\alpha$ is isomorphic to $M$ for
every $\alpha < \delta$ \then \, $M_\delta$ is isomorphic to $M$.
\end{definition}

\begin{remark}
\label{y.4} 
Concerning the following definition we 
shall use strongly limit in \ref{sl.7}(1), medium limit in 
\ref{sl.7}(2).
\end{remark}

\begin{definition}
\label{y.5}  
Let $\lambda$ be a cardinal $\ge |T|$.  
For parts 3) - 7) but not 8), for simplifying
the presentation we assume the axiom of global choice and $\bold F$ is
a class function; alternatively
restrict yourself to models with universe an ordinal $\in [\lambda,\lambda^+)$.

\noindent
1) For non-empty $\Theta \subseteq \{\mu:\aleph_0 \le \mu < \lambda$ and $\mu$
 is regular$\}$ and $M \in \text{ EC}_\lambda(T)$ we say that $M$ is a
$(\lambda,\Theta)$-superlimit \when \,: $M$ is universal and

if $\langle M_i:i \le \mu \rangle$ is $\prec$-increasing, 
$M_i \cong M$ for $i < \mu$ and  $\mu \in \Theta$ 

then $\cup \{M_i:i < \mu\} \cong M$.

\noindent
2) If $\Theta$ is a singleton, say $\Theta = \{\theta\}$, we may say that
$M$ is $(\lambda,\theta)$-superlimit.

\noindent
3) Let $S \subseteq \lambda^+$ be stationary.  A model 
$M \in \text{ EC}_\lambda(T)$
is called $S$-strongly limit or $(\lambda,S)$-strongly limit 
\when \, for some function:  $\bold F:\text{EC}_\lambda(T)
\rightarrow \text{ EC}_\lambda(T)$ we have:
\mn
\begin{enumerate}
\item[$(a)$]   for $N \in \text{ EC}_\lambda(T)$ 
we have $N \prec \bold F(N)$
\sn
\item[$(b)$]   if $\delta \in S$ is a limit ordinal
and $\langle M_i:i < \delta
\rangle$ is a $\prec$-increasing continuous sequence 
\footnote{no loss if we add $M_{i+1} \cong M$, so this simplifies the
demand on $\bold F$, i.e., only $\bold F(M')$ for $M' \cong M$ is required}
in EC$_\lambda(T)$ and 
$i < \delta \Rightarrow \bold F(M_{i+1}) \prec M_{i+2}$, 
\then \, $M \cong \cup\{M_i:i < \delta\}$.
\end{enumerate}
\mn
4) Let $S \subseteq \lambda^+$ be stationary.  $M \in 
\text{ EC}_\lambda(T)$ is called $S$-limit or $(\lambda,S)$-limit 
\underline{if} for some function $\bold F:\text{EC}_\lambda(T) \rightarrow
\text{ EC}_\lambda(T)$ we have:
\mn
\begin{enumerate}
\item[$(a)$]   for every $N \in \text{ EC}_\lambda(T)$ we have 
$N \prec \bold F(N)$ 
\sn
\item[$(b)$]    if $\langle M_i:i < \lambda^+ \rangle$ is
a $\prec$-increasing continuous sequence of members of 
EC$_\lambda(T)$ such that $\bold F(M_{i+1}) \prec M_{i+2}$ for $i <
\lambda^+$ \then \, 
for some closed unbounded \footnote{alternatively, we can use as a
parameter a filter on $\lambda^+$ extending the co-bounded filter} 
subset $C$ of $\lambda^+$,

\[
[\delta \in S \cap C \Rightarrow M_\delta \cong M].
\]
\end{enumerate}
\mn
5) We define\footnote{Note that $M$ is $(\lambda,S)$-strongly limit
iff $M$ is $(\{\lambda,\text{cf}(\delta):\delta \in S\})$-strongly limit.}
 ``$S$-weakly limit", ``$S$-medium limit" 
like ``$S$-limit", ``$S$-strongly limit" 
respectively by demanding that the domain of
$\bold F$ is the family of $\prec$-increasing continuous
sequence of members of EC$_\lambda(T)$ of length $< \lambda^+$ and
replacing ``$\bold F(M_{i+1}) \prec M_{i+2}$" by 
``$M_{i+1} \prec \bold F(\langle M_j:j \le i+1 \rangle) 
\prec M_{i+2}"$.  

\noindent
6) If $S = \lambda^+$ then we may omit $S$ (in parts (3), (4), (5)).  

\noindent
7) For non-empty $\Theta \subseteq \{\mu:\mu \le \lambda$ and $\mu$ is 
regular$\},M$ is $(\lambda,\Theta)$-strongly limit\footnote{in
\cite{Sh:88r} we consider:
we replace ``limit" by ``limit$^-$" if
$``\bold F(M_{i+1}) \prec M_{i+2}",``M_{i+1} \prec
\bold F(\langle M_j:j \le i+1 \rangle) \prec M_{i+2}"$ are replaced 
by $``\bold F(M_i) \prec M_{i+1}",``M_i \prec 
\bold F(\langle M_j:j \le i \rangle) \prec M_{i+1}"$ respectively.
But $(\text{EC}(T),\prec)$ has amalgamation.} 
\underline{if} $M$ is 
$\{\delta < \lambda^+:\text{cf}(\delta) \in \Theta\}$-strongly limit.  
Similarly for the other notions.  If we do not write $\lambda$
we mean $\lambda = \|M\|$.  

\noindent
8) We say that $M \in K_\lambda$ is invariantly strong limit \when \, in
part (3), $\bold F$ is just a subset of $\{(M,N)/\cong:M
\prec N$ are from EC$_\lambda(T)\}$ and in clause (b) of part (3)
we replace ``$\bold F(M_{i+1}) \prec M_{i+2}"$ by ``$(\exists
N)(M_{i+1} \prec N \prec M_{i+2} \wedge ((M,N)/\cong) \in 
\bold F)$".  But abusing notation we still 
write $N = \bold F(M)$ instead $((M,N)/ \cong) \in \bold F$.  
Similarly with the other notions, so we use the isomorphism type of
$\bar M \char 94 \langle N\rangle$ for ``weakly limit" and ``medium
limit".

\noindent
9) In the definitions above we may say ``$\bold F$ witness $M$ is ..."  
\end{definition}

\begin{observation}
\label{y.5c}  
1) Assume $\bold F_1,\bold F_2$ are as above and
$\bold F_1(N) \prec \bold F_2(N)$ (or $\bold F_1(\bar N) \prec \bold
F_2(\bar N))$ whenever defined.  If $\bold F_1$ is a witness then 
so is $\bold F_2$.

\noindent
2) All versions of limit models implies being a universal model in 
{\rm EC}$_\lambda(T)$. 
\end{observation}
\bigskip

\noindent
3) \underline{The Obvious implications diagram}:
\label{y.6}  
For non-empty $\Theta \subseteq \{\theta:\theta$ is regular $\le
\lambda\}$ and stationary $S_1 \subseteq \{\delta <
\lambda^+:\text{cf}(\delta) \in \Theta\}$:

\[
\text{superlimit } = (\lambda,\{\mu:\mu \le \lambda \text{
regular}\})\text{-superlimit}
\]
\centerline {$\downarrow$}
\[
(\lambda,\Theta)\text{-superlimit}
\]
\centerline {$\downarrow$}
\[
S_1 \text{-strongly limit}
\]
\centerline {$\downarrow$ \hskip25pt $\downarrow$}

\hskip75pt $S_1 \text{-medium limit}, \qquad \qquad \qquad \quad 
S_1\text{-limit}$

\centerline {$\downarrow$ \hskip25pt $\downarrow$}
\[
S_1 \text{-weakly limit}.
\]

\begin{lemma}
\label{y.7}  
Let $T$ be a first order complete theory.

\noindent
1) If $\lambda$ is regular, $M$ a saturated model of $T$ of cardinality
$\lambda$, \then \, $M$ is $(\lambda,\lambda)$-superlimit. 

\noindent
2) If $T$ is stable, and $M$ is a saturated model of $T$ of cardinality
$\lambda \ge \aleph_1 + |T|$ and $\Theta = \{\mu:\kappa(T) \le \mu \le
\lambda$ and $\mu$ is regular$\})$, \then \, $M$ is 
$(\lambda,\Theta)$-superlimit (on
$\kappa(T)$-see \cite[III,\S3]{Sh:c}).   

\noindent
3) If $T$ is stable in $\lambda$ and $\kappa = \text{\rm cf}(\kappa) \le
\lambda$ \then \, $T$ has an invariantly strongly 
$(\lambda,\kappa)$-limit model.
\end{lemma}

\begin{remark}
\label{y8}
Concerning \ref{y.7}(2), note that by \cite{Sh:c} if
$\lambda$ is singular or just $\lambda < \lambda^{< \lambda}$
and $T$ has a saturated model of cardinality
$\lambda$ \then \, $T$ is stable (even stable in $\lambda$) 
and $\cf(\lambda) \ge \kappa(T))$.
\end{remark}

\begin{PROOF}{\ref{y.7}}
1) Let $M_i$ be a $\lambda$-saturated model of $T$ of cardinality 
$\lambda$ for $i < \lambda$
and $\langle M_i:i < \lambda\rangle$ is $\prec$-increasing and
$M_\lambda = \bigcup\limits_{i < \lambda} M_i$.  Now
for every $A \subseteq M_\lambda$ of cardinality $< \lambda$
there is $i < \lambda$ such that $A \subseteq M_i$ hence every $p \in
\bold S(A,M_\lambda)$ is realized in $M_i$ hence in $M_\lambda$; so
clearly $M_\lambda$ is $\lambda$-saturated.  Remembering the uniqueness of a
$\lambda$-saturated model of $T$ of cardinality $\lambda$ we finish. 

\noindent
2) Use \cite[III,3.11]{Sh:c}: if $M_i$ is a $\lambda$-saturated model
of $T,\langle M_i:i < \delta \rangle$ increasing cf$(\delta) \ge
\kappa(T)$ \then \, $\bigcup\limits_{i < \delta} M_i$ is $\lambda$-saturated. 

\noindent
3)  Let $\bold K_{\lambda,\kappa} = \{\bar M:\bar M = \langle M_i:i \le
    \kappa\rangle$ is $\prec$-increasing continuous, $M_i \in \text{
    EC}_\lambda(T)$ and $(M_{i+2},c)_{c \in M_{i+1}}$ is saturated for
    every $i < \kappa\}$.  Clearly $\bar M,\bar N \in \bold
    K_{\lambda,\kappa} \Rightarrow M_\kappa \cong N_\kappa$.  Also for
    every $M \in \text{ EC}_\lambda(T)$ there is $N$ such that $M
    \prec N$ and $(N,c)_{c \in M}$ is saturated, as also Th$((M,c)_{c
    \in M})$ is stable in $\lambda$; so there is an invariant
$\bold F:\text{EC}_\lambda(T) \rightarrow \text{ EC}_\lambda(T)$ such
    that $M \prec \bold F(M)$ and $(\bold F(M),c)_{c \in M}$ is
    saturated; such $\bold F$ witness the desired conclusion.
\end{PROOF}

\begin{definition}
\label{y15}
0) For regular $\kappa < \lambda$ let $S^\lambda_\theta = \{\delta <
   \lambda:\cf(\delta) = \lambda\}$.

\noindent
1) For a regular uncountable cardinal $\lambda$ let 
$\check I[\lambda] = \{S \subseteq \lambda$: some
pair $(E,\bar a)$ witnesses $S \in \check I[\lambda]$, see below$\}$. 

\noindent
2) We say that $(E,\bar u)$ is a witness for $S \in \check I[\lambda]$
\Iff :
\mn
\begin{enumerate}
\item[$(a)$]   $E$ is a club of the regular cardinal $\lambda$
\sn
\item[$(b)$]   $\bar u = \langle u_\alpha:\alpha < \lambda
\rangle,u_\alpha \subseteq \alpha$ and $\beta \in u_\alpha \Rightarrow
u_\beta = \beta \cap u_\alpha$
\sn
\item[$(c)$]   for every $\delta \in E \cap S,u_\delta$ is an
unbounded subset of $\delta$ of order-type cf$(\delta)$ (and $\delta$ is
a limit ordinal).
\end{enumerate}
\end{definition}

By \cite[\S1]{Sh:420}
\begin{claim}
\label{y16}
If $\kappa^+ < \lambda$ and $\kappa,\lambda$ are regular then some
stationary $S \subseteq \{\delta < \lambda:\text{cf}(\delta) =
\kappa\}$ belongs to $\check I[\lambda]$.
\end{claim}

\noindent
By \cite{Sh:108}
\begin{claim}
\label{y18}
If $\lambda = \mu^+,\theta = \cf(\theta) \le \cf(\mu)$ 
and $\alpha < \mu \Rightarrow
|\alpha|^{< \theta} \le \mu$ \then \, $S^\lambda_\theta \in \check I[\lambda]$.
\end{claim}
\newpage

\section {On superstable not $\aleph_0$-stable $T$} \label{Onsuperstable}

We first note that superstable $T$ tend to have superlimit models.
\begin{claim}
\label{nlm.0.0}  
Assume $T$ is superstable and $\lambda
\ge |T| + 2^{\aleph_0}$.  Then $T$ has a superlimit model of
cardinality $\lambda$ \underline{iff} $T$ has a saturated model of
cardinality $\lambda$ \underline{iff} $T$ has a universal model of cardinality
$\lambda$ iff $\lambda \ge |D(T)|$. 
\end{claim}

\begin{PROOF}{\ref{nlm.0.0}}
By \cite[III,\S5]{Sh:c} we know that $T$ is stable in
$\lambda$ iff $\lambda \ge |D(T)|$.  Now if $|T| \le
\lambda < |D(T)|$ trivially there is no universal model of $T$ of
cardinality $\lambda$ hence no saturated model and no superlimit
model, etc., recalling \ref{y.5c}(2).  
If $\lambda \ge |D(T)|$, then $T$ is stable in $\lambda$
hence has a saturated model of cardinality $\lambda$ by \cite[III]{Sh:c}
(hence universal)
and the class of $\lambda$-saturated models of $T$ is closed under
increasing elementary chains by \cite[III]{Sh:c} so we are done.
\end{PROOF}

\noindent
The following are the prototypical theories which we shall consider.
\begin{definition}
\label{0.0}  
1) $T_0 = \text{ Th}({}^\omega 2,E^0_n)_{n < \omega}$ 
when $\eta E^0_n \nu \Leftrightarrow \eta \restriction n = \nu \restriction n$.

\noindent
2) $T_1 = \text{ Th}({}^\omega(\omega_1),E^1_n)_{n < \omega}$ where 
$\eta E^1_n \nu \Leftrightarrow \eta \restriction n = \nu \restriction n$.

\noindent
3) $T_2 = \text{ Th}(\Bbb R,<)$.
\end{definition}

\noindent
Recall
\begin{observation}
\label{0.1.1}  
0) $T_\ell$ is a countable complete first order theory for
   $\ell=0,1,2$.

\noindent
1) $T_0$ is superstable not $\aleph_0$-stable.

\noindent
2) $T_1$ is strictly stable, that is, stable not superstable.

\noindent
3) $T_2$ is unstable.

\noindent
4) $T_\ell$ has elimination of quantifiers for $\ell=0,1,2$.
\end{observation}

\begin{claim}
\label{0.1.2}  
It is consistent with {\rm ZFC} that
$\aleph_1 < 2^{\aleph_0}$ and some $M \in \text{\rm EC}_{\aleph_1}
(T_0)$ is a superlimit model.
\end{claim}

\begin{PROOF}{\ref{0.1.2}}
By \cite{Sh:100}, for notational simplicity we start with $\bold V =
\bold L$.  

So $T_0$ is defined in \ref{0.0}(1) and it is the $T$ from 
Theorem \cite[1.1]{Sh:100} and let $S$ be the set of $\eta \in ({}^\omega
2)^{\bold L}$.  We define $T'$ (called $T_1$ there) as the following theory:
\mn
\begin{enumerate}
\item[$\circledast_1$]  $(i) \quad T_0$, or just for each $n$ the sentence
saying $E_n$ is an equivalence 

\hskip25pt relation with $2^n$ equivalence classes, 
each $E_n$ equivalence class

\hskip25pt  divided to two by $E_{n+1},E_{n+1}$ refine $E_n,E_0$ is trivial
\sn
\item[${{}}$]   $(ii) \quad$ the sentences saying that
\sn
\begin{enumerate}
\item[${{}}$]   $(\alpha) \quad$ for every $x$, the function $z
\mapsto F(x,z)$ is one-to-one and 
\sn
\item[${{}}$]   $(\beta) \quad x E_n(F(x,z))$ for each $n < \omega$
\end{enumerate}
\item[${{}}$]   $(iii) \quad E_n(c_\eta,c_\nu)^{\text{if}(\eta
\restriction n=\nu \restriction n)}$ for $\eta,\nu \in S$.
\end{enumerate}
\mn
In \cite{Sh:100} it is proved that in some forcing\footnote{We can
replace $\bold L$ by any $\bold V_0$ which satisfies $2^{\aleph_0} =
\aleph_1,2^{\aleph_1} = \aleph_2$.}
extension $\bold L^{\bbP}$ of $\bold L$, $\bbP$ an
$\aleph_2$-c.c. proper forcing of cardinality $\aleph_2$, 
in $\bold V = \bold L^{\bbP}$, the class
PC$(T',T_0) = \{M \restriction \tau_{T_0}:M$ 
is a $\tau$-model of $T'\}$ is categorical in $\aleph_1$.

However, letting $M^*$ be any model from PC$(T',T_0)$ of cardinality
$\aleph_1$, it is easy to see that (in $\bold V = \bold L^{\bbP}$):
\mn
\begin{enumerate}
\item[$\circledast_2$]   the following conditions on $M$ are
equivalent
\begin{enumerate}
\item[$(a)$]   $M$ is isomorphic to $M^*$
\sn
\item[$(b)$]    $M \in \text{\rm PC}(T',T_0)$ 
\sn
\item[$(c)$]   $(\alpha) \quad M$ is a model of $T_0$ of 
cardinality $\aleph_1$
\sn
\item[${{}}$]   $(\beta) \quad M^*$ can be elementarily embedded into $M$
\sn
\item[${{}}$]    $(\gamma) \quad$ for every 
$a \in M$ the set $\cap\{a/E^M_n:n < \omega\}$ has cardinality $\aleph_1$.
\end{enumerate}
\mn
But
\mn
\item[$\circledast_3$]   every model $M_1$ of $T$ of cardinality
$\le \aleph_1$ has a proper elementary extension to a model satisfying
(c), i.e., $(\alpha),(\beta),(\gamma)$ of $\circledast_2$ above
\sn
\item[$\circledast_4$]   if $\langle M_\alpha:\alpha < \delta
\rangle$ is an increasing chain of models satisfying (c) of
$\circledast_2$ and $\delta < \omega_2$ then also
$\cup\{M_\alpha:\alpha < \delta\}$ does.
\end{enumerate}
\mn
Together we are done.  
\end{PROOF}

\noindent
Naturally we ask
\begin{question}
\label{0.1.3}  
What occurs to $T_0$ for $\lambda > \aleph_1$ but $\lambda < 2^{\aleph_0}$?
\end{question}

\begin{question}
\label{0.1}  
Does the theory $T_2$ of linear order consistently have an
$(\aleph_1,\aleph_0)$-superlimit? (or only strongly limit?) but see \S3.
\end{question}

\begin{question}
\label{0.2}  
What is the answer for $T$ when $T$ is 
countable superstable not $\aleph_0$-stable and $D(T)$ countable for
$\aleph_1 < 2^{\aleph_0}$ for $\aleph_2 < 2^{\aleph_0}$?

So by the above for some such $T$, in some universe, for $\aleph_1$ the
answer is yes, there is a superlimit.
\end{question}
\newpage

\section {A strictly stable consistent example} \label{Astrictly}

We now look at models of $T_1$ (redefined below) in cardinality
$\aleph_1$; recall
\begin{definition}
\label{sl.1}
$T_1 = \text{ Th}({}^\omega(\omega_1),E_n)_{n < \omega}$ where $E_n =
\{(\eta,\nu):\eta,\nu \in {}^\omega(\omega_1)$ and $\eta \restriction
n = \nu \restriction n\}$.
\end{definition}

\begin{remark}
\label{sl.1.3}  
\mn
\begin{enumerate}
\item[$(a)$]   Note that $T_1$ has elimination of quantifiers.
\sn
\item[$(b)$]   If $\lambda = \Sigma\{\lambda_n:n < \omega\}$
and $\lambda_n = \lambda_n^{\aleph_0}$, \then \, $T_1$ has a
$(\lambda,\aleph_0)$-superlimit model in $\lambda$ (see \ref{sl.11}).
\end{enumerate}
\end{remark}

\begin{defclaim}
\label{sl.1y}  
1) Any model of $T_1$ of cardinality $\lambda$ is isomorphic to
$M_{A,h} := (\{(\eta,\varepsilon):\eta \in A,
\varepsilon < h(\eta))\},E_n)_{n < \omega}$ 
for some $A \subseteq {}^\omega \lambda$ and 
$h:{}^\omega \lambda \rightarrow (\Car \cap \lambda^+) \backslash
\{0\}$ where $(\eta_1,\varepsilon_1) E_n(\eta_2,\varepsilon_2) \Leftrightarrow 
\eta_1 \restriction n = \eta_2 \restriction n$, pedantically we should write
$E^{M_{A,h}}_n = E_n \rest |M_{A,n}|$.

\noindent
2) We write $M_A$ for $M_{A,h}$ when $A$ is as above and $h:A
\rightarrow \{|A|\}$, so constantly $|A|$ when $A$ is infinite.

\noindent
3) For $A \subseteq {}^\omega \lambda$ and $h$ as above the model
   $M_{A,h}$ is a model of $T_1$ \Iff \, $A$ is non-empty and $(\forall
   \eta \in A)(\forall n < \omega)(\exists^{\aleph_0} \nu \in A)(\nu
   \rest n =\eta \rest n \wedge \nu(n) \ne \eta(n))$.

\noindent
4) Above $M_{A,h}$ has cardinality $\lambda$ iff $\Sigma\{h(\eta):\eta
   \in A\} = \lambda$.
\end{defclaim}

\begin{definition}
\label{sl.2}  
1) We say that $A$ is a $(T_1,\lambda)$-witness \when
\mn
\begin{enumerate}
\item[$(a)$]   $A \subseteq {}^\omega \lambda$ has cardinality $\lambda$
\sn
\item[$(b)$]   if $B_1,B_2 \subseteq {}^\omega \lambda$ are
$(T_1,A)$-big (see below) of cardinality $\lambda$ \then \,
$(B_1 \cup {}^{\omega >}\lambda,\triangleleft)$ is isomorphic to
$(B_2 \cup {}^{\omega >}\lambda,\triangleleft)$.
\end{enumerate}
\mn
2) A set $B \subseteq {}^\omega \lambda$ is called $(T_1,A)$-big \when
\, it is $(\lambda,\lambda)-(T_1,A)$-big; see below.

\noindent
3) $B$ is $(\mu,\lambda)-(T_1,A)$-big means: $B \subseteq {}^\omega \lambda,
|B| = |A| = \mu$ and
 for every $\eta \in {}^{\omega >}\lambda$ 
there is an isomorphism $f$ from $({}^{\omega \ge}\lambda,\triangleleft)$ onto
$(\{\eta \char 94 \nu:\nu \in{}^{\omega \ge}\lambda\},\triangleleft)$
mapping $A$ into $\{\nu:\eta \char 94 \nu \in B\}$.

\noindent
4) $A \subseteq {}^\omega(\omega_1)$ is $\aleph_1$-suitable \when \,:
\mn
\begin{enumerate}
\item[$(a)$]   $|A| = \aleph_1$
\sn
\item[$(b)$]   for a club of $\delta < \omega_1,A \cap{}^\omega
\delta$ is everywhere not meagre in the space ${}^\omega \delta$, i.e.,
for every $\eta \in {}^{\omega >}\delta$ the set $\{\nu \in A \cap
{}^\omega \delta:\eta \triangleleft \nu\}$ is a non-meagre subset of
${}^\omega \delta$ (that is what really is used in \cite{Sh:100}).
\end{enumerate}
\end{definition}

\begin{claim}
\label{st.3}  
It is consistent with ZFC that $2^{\aleph_0} > 
\aleph_1 +$ there is a $(T_1,\aleph_1)$-witness;
moreover every $\aleph_1$-suitable set is a $(T_1,\aleph_1)$-witness.
\end{claim}

\begin{PROOF}{\ref{st.3}}
By \cite[\S2]{Sh:100}.
\end{PROOF}

\begin{remark}
\label{st.3a}  
The witness does not give rise to an
$(\aleph_1,\aleph_0)$-limit model, as for the union of any ``fast enough"
$\prec$-increasing $\omega$-chain of members of EC$_{\aleph_1}(T_1)$,
the relevant sets are meagre.
\end{remark}

\begin{definition}
\label{sl.4}  
Let $A$ be a $(T_1,\lambda)$-witness.  We define $K^1_{T_1,A}$ as the family of
$M = (|M|,<^M,P^M_\alpha)_{\alpha \le \omega}$ such that:
\mn
\begin{enumerate}
\item[$(\alpha)$]    $(|M|,<^M)$ is a tree with $(\omega +1)$ levels 
\sn
\item[$(\beta)$]   $P^M_\alpha$ is the $\alpha$-th level; let $P^M_{<
\omega} = \cup\{P^M_n:n < \omega\}$
\sn
\item[$(\gamma)$]   $M$ is isomorphic to $M^1_B$ for some 
$B \subseteq {}^\omega \lambda$ of cardinality $\lambda$ where $M^1_B$
is defined by $|M^1_\beta| = ({}^{\omega >} \lambda)
\cup B,P^{M^1_B}_n = {}^n \lambda,P^{M^1_B}_\omega = B$ 
and $<^{M^1_B} = \triangleleft  \rest |M^1_B|$, i.e., being an initial segment
\sn
\item[$(\delta)$]   moreover $B$ is such that some $f$ satisfies:
\begin{enumerate}
\item[$\circledast$]   $(a) \quad f:{}^{\omega >}\lambda
\rightarrow \omega$ and $f(<>) = 0$ for simplicity
\sn
\item[${{}}$]    $(b) \quad \eta \trianglelefteq \nu \in {}^{\omega
>}\lambda \Rightarrow f(\eta) \le f(\nu)$
\sn
\item[${{}}$]    $(c) \quad$ if $\eta \in B$ then $\langle f(\eta
\restriction n):n < \omega \rangle$ is eventually constant
\sn
\item[${{}}$]   $(d) \quad$ if $\eta \in {}^{\omega >}\lambda$
then $\{\nu \in {}^\omega\lambda:\eta {}^\frown \nu \in B$ and $m <
\omega \Rightarrow$

\hskip30pt $f(\eta {}^\frown (\nu \restriction m))= 
f(\eta)\}$ is $(T_1,A)$-big
\sn
\item[${{}}$]   $(e) \quad$ for $\eta \in {}^{\omega >}\lambda$
and $n \in [f(\eta),\omega)$ for $\lambda$ ordinals $\alpha <
\lambda$, we have 

\hskip30pt  $f(\eta {}^\frown \langle \alpha \rangle)=n$.
\end{enumerate}
\end{enumerate}
\end{definition}

\begin{claim}
\label{sl.5}
[The Global Axiom of Choice]   
If $A$ is a $(T_1,\aleph_1)$-witness \then
\mn
\begin{enumerate}
\item[$(a)$]   $K^1_{T_1,A} \ne \emptyset$
\sn
\item[$(b)$]   any two members of $K^1_{T_1,A}$ are isomorphic
\sn
\item[$(c)$]   there is a function $\bold F$ from $K^1_{T_1,A}$ to
itself (up to isomorphism, i.e., $(M,\bold F(M))$ is defined only up
to isomorphism) satisfying $M \subseteq \bold F(M)$ 
such that $K^1_{T_1,A}$ is closed
under increasing unions of sequence $\langle M_n:n < \omega
\rangle$ such that $\bold F(M_n) \subseteq M_{n+1}$.
\end{enumerate}
\end{claim}

\begin{PROOF}{\ref{sl.5}}
\smallskip

\noindent
\underline{Clause (a)}:  Trivial.
\smallskip

\noindent
\underline{Clause (b)}:  By the definition of ``$A$ is a
$(T_1,\aleph_1)$-witness" and of $K^1_{T_1,A}$.
\smallskip

\noindent
\underline{Clause (c)}:
\smallskip

We choose $\bold F$ such that
\mn
\begin{enumerate}
\item[$\circledast$]   if $M \in K^1_{A,T_1}$ then $M \subseteq
\bold F(M) \in K^1_{A,T_1}$ and for every $k < \omega$ and $a \in P^M_k$, the
set $\{b \in P^{\bold F(M)}_{k+1}:a <_{\bold F(M)} b$ and $b \notin M\}$ has
cardinality $\aleph_1$.
\end{enumerate}
\mn
Assume $M = \cup\{M_n:n < \omega\}$ where $\langle M_n:n <\omega\rangle$ is
$\subseteq$-increasing$\}$, $M_n \in K^1_{A,T_1},\bold F(M_n)
\subseteq M_{n+1}$.  Clearly $M$ is as required in the beginning of
Definition \ref{sl.4}, that is, satisfies clauses
$(\alpha),(\beta),(\gamma)$ there.  To prove clause $(\delta)$, 
we define $f:P^M_{<\omega} \rightarrow
\omega$ by $f(a) = \text{ Min}\{n:a \in M_n\}$.  Pendantically, $\bold
F$ is defined only up to isomorphism.

\noindent
So we are done.   
\end{PROOF}

\begin{claim}
\label{sl.5.1}
[The Global Axiom of Choice]

If $A$ is a $(T_1,\lambda)$-witness \then \,
\mn
\begin{enumerate}
\item[$(a)$]   $K^1_{T_1,A} \ne \emptyset$
\sn
\item[$(b)$]   any two members of $K^1_{T_1,A}$ are isomorphic
\sn
\item[$(c)$]   if $M_n \in K^1_{T_1,A}$ and $n < \omega \Rightarrow
M_n \subseteq M_{n+1}$ \then \, $M := \cup\{M_n:n < \omega\} \in
K^1_{T_1,A}$.
\end{enumerate}
\end{claim}

\begin{remark}
If we omit clause (b), we can weaken the demand on the set $A$.
\end{remark}

\begin{PROOF}{\ref{sl.5.1}}
Assume $M = \cup\{M_n:n < \omega\},M_n
\subseteq M_{n+1},M_n \in K^1_{T_1,A}$ and $f_n$ witnesses $M_n \in
K^1_{T_1,A}$.  Clearly $M$ satisfies clauses
$(\alpha),(\beta),(\gamma)$ from Definition \ref{sl.4}, we just 
have to find a witness $f$ as in clause $(\delta)$ there.

For each $a \in M$ let $n(a) = \text{ Min}\{n:a \in M_n\}$, clearly if
$M \models ``a < b < c"$ then $n(a) \le n(b)$ and $n(a) = n(c)
\Rightarrow n(a) = n(b)$.  Let $g_n:M \rightarrow M$ be defined by: 
$g_n(a)=b$ iff $b \le^M a,b \in M_n$ and $b$ is $\le^M$-maximal under those
restrictions; clearly it is well defined.  Now we define $f'_n:M_n
\rightarrow \omega$ by induction on $n < \omega$ such that $m<n
\Rightarrow f'_m \subseteq f'_n$, as follows.

If $n=0$ let $f'_n = f_n$.

If $n=m+1$ and $a \in M_n$ we let $f'_n(a)$ be $f'_m(a)$ if $a \in
M_m$ and be $(f_n(a)-f_n(g_m(a))) + f'_m(g_m(a))+1$ if $a \in M_n
\backslash M_m$.  Clearly $f := \cup\{f'_n:n < \omega\}$ is a function
from $M$ to $\omega,a \le^M b \Rightarrow f(a) \le f(b)$, and for any
$a \in M$ the set $\{b \in M:a \le^M b$ and $f(b) = f(a)\}$ is equal
to $\{b \in M_{n(a)}:f_{n(a)}(a) = f_{n(a)}(b)$ and $a \le^M b\}$.  

\noindent
So we are done.
\end{PROOF}

\begin{definition}
\label{sl.5.7}  Let $A$ be a $(T_1,\lambda)$-witness.  
We define $K^2_{T_1,A}$ as in Definition \ref{sl.4} but $f$ is constantly zero.
\end{definition}

\begin{claim}
\label{sl.6} 
[The Global Axiom of Choice] 
 If $A$ is a $(T_1,\aleph_1)$-witness \then \,
\mn
\begin{enumerate}
\item[$(a)$]  $K^2_{T_1,A} \ne \emptyset$
\sn
\item[$(b)$]   any two members of $K^2_{T_1,A}$ are isomorphic
\sn
\item[$(c)$]   there is a function $\bold F$ from 
$\cup\{{}^{\alpha +2}(K^2_{T_0,A}): \alpha < \omega_1\}$ to
$K^2_{T_1,A}$ which satisfies:
\begin{enumerate}
\item[$\boxtimes$]   $(\alpha) \quad$ if $\bar M = \langle M_i:i \le
\alpha + 1 \rangle$ is an $\prec$-increasing sequence of models

\hskip30pt  of $T$ then $M_{\alpha +1} \subseteq 
\bold F(\bar M) \in K^2_{T_1,A}$
\sn
\item[${{}}$]   $(\beta) \quad$  the union of
any increasing $\omega_1$-sequence $\bar M = \langle M_\alpha:\alpha <
\omega_1 \rangle$

\hskip30pt  of members of $K^2_{T_1,A}$ belongs to $K^2_{T_1,A}$
when 

\hskip30pt  $\omega_1 = \sup\{\alpha:\bold F(\bar M \restriction (\alpha+2))
\subseteq M_{\alpha +2})$ and is a well defined

\hskip30pt   embedding of $M_\alpha$ into $M_{\alpha +2}\}$.
\end{enumerate}
\end{enumerate}
\end{claim}

\begin{remark}
Instead of the global axiom of choice, we can
restrict the models to have universe a subset of $\lambda^+$ (or just
a set of ordinals).
\end{remark}

\begin{PROOF}{\ref{sl.6}}  
\smallskip

\noindent
\underline{Clause (a)}:  Easy.
\smallskip

\noindent
\underline{Clause (b)}:  By the definition.
\smallskip

\noindent
\underline{Clause (c)}:  Let 
$\langle {\cU}_\varepsilon:\varepsilon < \omega_1
\rangle$ be an increasing sequence of subsets of $\omega_1$
with union $\omega_1$ such that $\varepsilon < \omega_1 \Rightarrow
|{\cU}_\varepsilon \backslash \bigcup\limits_{\zeta < \varepsilon} 
{\cU}_\zeta| = \aleph_1$.  Let $M^* \in K^2_{T_1,A}$
be such that ${}^{\omega >}(\omega_1) \subseteq |M^*| \subseteq
{}^{\omega \ge}(\omega_1)$ and $M^*_\varepsilon =: M^* \restriction
{}^{\omega \ge}({\cU}_\varepsilon)$ belongs to $K^2_{T_1,A}$ for
every $\varepsilon < \omega_1$.

We choose a pair $(\bold F,\bold f)$ of functions with domain $\{\bar
M:\bar M$ an increasing sequence of members of $K^2_{T_1,A}$
of length $< \omega_1\}$ such that:
\mn
\begin{enumerate}
\item[$(\alpha)$]   $\bold F(\bar M)$ is an extension of
$\cup\{M_i:i < \ell g(\bar M)\}$ from $K^2_{T_1,A}$
\sn
\item[$(\beta)$]   $\bold f(\bar M)$ is an embedding from
$M^*_{\ell g(\bar M)}$ into $\bold F(\bar M)$
\sn
\item[$(\gamma)$]   if $\bar M^\ell = \langle M_\alpha:\alpha <
\alpha_\ell \rangle$ for $\ell=1,2$ and 
$\alpha_1 < \alpha_2,\bar M^1 = \bar M^2 \restriction
\alpha_1$ and $\bold F(\bar M^1) \subseteq M_{\alpha_1}$
then $\bold f(\bar M^1) \subseteq \bold f(\bar M^2)$
\sn
\item[$(\delta)$]   if $a \in \bold F(\bar M)$ and $n < \omega$ then
for some $b \in M^*_{\ell g(\bar M)}$ we have $\bold F(M) \models a E_n
(\bold f(\bar M)(b))$.  
\end{enumerate}
\mn
Now check.
\end{PROOF}

\begin{conclusion}
\label{sl.7} 
Assume there is a $(T_1,\aleph_1)$-witness (see Definition \ref{sl.2})
for the first-order complete theory $T_1$ from \ref{sl.1}:

\noindent
1) $T_1$ has an $(\aleph_1,\aleph_0)$-strongly limit model.

\noindent
2)  $T_1$ has an $(\aleph_1,\aleph_1)$-medium limit model.

\noindent
3) $T_1$ has a $(\aleph_1,\aleph_0)$-superlimit model.
\end{conclusion}

\begin{PROOF}{\ref{sl.7}}
1)  By \ref{sl.5} the reduction of problems on
$(\text{EC}(T_1),\prec)$ to $K^1_{T_1,A}$ (which is easy) is
exactly as in \cite{Sh:100}.

\noindent
2) By \ref{sl.6}.

\noindent
3) Like part (1) using claim \ref{sl.5.1}.  
\end{PROOF}

\begin{claim}
\label{sl.11}
If $\lambda = \Sigma\{\lambda_n:n < \omega\}$
and $\lambda_n = \lambda_n^{\aleph_0}$, \then \, $T_1$ has a
$(\lambda,\aleph_0)$-superlimit model in $\lambda$.
\end{claim}

\begin{PROOF}{\ref{sl.11}}
Let $M_n$ be the model $M_{A_n,h_n}$ where $A_n =
{}^\omega(\lambda_n)$ and $h_n:A_n \rightarrow \lambda^+_n$ is
constantly $\lambda_n$.

Clearly
\mn
\begin{enumerate}
\item[$(*)_1$]   $M_n$ is a saturated model of $T_1$ of cardinality
$\lambda_n$
\sn
\item[$(*)_2$]  $M_n \prec M_{n+1}$
\sn
\item[$(*)_3$]  $M_\omega = \cup\{M_n:n <\omega\}$ is a special model
of $T_1$ of cardinality $\lambda$.
\end{enumerate}
\mn
The main point:
\mn
\begin{enumerate}
\item[$(*)_4$]  $M_\omega$ is $(\lambda,\aleph_0)$-superlimit model of
$T_1$.
\end{enumerate}
\mn
[Why?  Toward this assume
\mn
\begin{enumerate}
\item[$(a)$]  $N_n$ is isomorphic to $M_\omega$ say $f_n:M_\omega
\rightarrow N_n$ is such isomorphic
\sn
\item[$(b)$]  $N_n \prec N_{n+1}$ for $n < \omega$.
\end{enumerate}
\mn
Let $N_\omega = \cup\{N_n:n < \omega\}$ and we should prove $N_\omega
\cong M_\omega$, so just $N_\omega$ is a special model of $T_1$ of
cardinality $\lambda$ suffice.

Let $N'_n = N_\omega \rest (\cup\{f_n(M_k):k \le n\})$.  Easily $N'_n
\prec N'_{n+1} \prec N_\omega$ and $\cup\{N'_n:n < \omega\} =
N_{\omega_*}$ and $\|N'_n\| = \lambda_n$.  So it suffices to prove
that $N'_n$ is saturated and by direct inspection shows this.
\end{PROOF}
\newpage

\section {On non-existence of limit models} \label{On non}

Naturally we assume that non-existence of superlimit models for unstable $T$ is
easier to prove.  For other versions we need to look more.  
We first show that for $\lambda \ge
|T| + \aleph_1$, if $T$ is unstable then it does not have a superlimit
model of cardinality $\lambda$ and if $T$ is 
unsuperstable, we show this for ``most" cardinals $\lambda$.  On
``$\Phi$ proper for $K_{\text{or}}$ or $K^\omega_{\text{tr}}$", see
\cite[VII]{Sh:c} or \cite{Sh:E59} or hopefully some day in
\cite[III]{Sh:e}.  We assume some knowledge on stability.
\begin{claim}
\label{nlm.0.1} 
1) If $T$ is unstable, $\lambda \ge |T| + \aleph_1$,
\then \, $T$ has no superlimit model of cardinality $\lambda$.

\noindent
2) If $T$ is stable not superstable and $\lambda \ge |T| + \beth_\omega$ or
$\lambda = \lambda^{\aleph_0} \ge |T|$ \then \,
$T$ has no superlimit model of cardinality $\lambda$.
\end{claim}

\begin{remark} 
1) We assume some knowledge on EM models for linear
orders $I$ and members of $K^\omega_{\text{tr}}$ as index models, see,
e.g. \cite[VII]{Sh:c}.

\noindent
2) We use the following definition in the proof, as well as 
a result from \cite{Sh:460} or \cite{Sh:829}. 
\end{remark}

\begin{definition}
\label{c8}
For cardinals $\lambda > \kappa$ let $\lambda^{[\kappa]}$ be the
minimal $\mu$ such that for some, equivalently
for every set $A$ of cardinality $\lambda$ there is ${\cP}_A \subseteq
[A]^{\le \kappa} = \{B \subseteq A:|B| \le \kappa\}$ of cardinality
$\lambda$ such that any $B \in [\lambda]^{\le \kappa}$ is the union of
$< \kappa$ members of ${\cP}_A$.
\end{definition}

\begin{PROOF}{\ref{nlm.0.1}}  
1) Towards a contradiction assume $M^*$ is a superlimit
model of $T$ of cardinality $\lambda$.  As $T$ is unstable we can find
$m,\varphi(\bar x,\bar y)$ such that
\mn
\begin{enumerate}
\item[$(*)$]  $\varphi(\bar x,\bar y) \in \bbL_{\tau(T)}$ 
linearly orders some infinite $\bold I \subseteq {}^m M,M \models T$
so $\ell g(\bar x) = \ell g(\bar y) = m$.
\end{enumerate}
\mn
We can find a $\Phi$ which is proper for linear orders (see
\cite[VII]{Sh:c}) and $F_\ell(\ell < m)$ such that 
$F_\ell \in \tau_\Phi \backslash \tau_T$ is a unary 
function symbol for $\ell < m,
\tau_T \subseteq \tau(\Phi)$ and for
every linear order $I$, EM$(I,\Phi)$ has Skolem functions and its
$\tau_T$-reduct EM$_{\tau(T)}(I,\Phi)$ is a 
model of $T$ of cardinality $|T| + |I|$ and 
$\tau(\Phi)$ is of cardinality $|T| + \aleph_0$ and
$\langle a_s:s \in I\rangle$ is the Skeleton of EM$(I,\Phi)$, that is,
it is an indiscernible sequence in EM$(I,\Phi)$ and EM$(I,\Phi)$ is
the Skolem hull of $\{a_s:s \in I\}$, and letting $\bar a_s = \langle
F_\ell(a_s):\ell < m\rangle$ in EM$(I,\Phi)$ we have
EM$_{\tau(T)}(I,\Phi) \models \varphi 
[\bar a_s,\bar a_t]^{\text{if}(s<t)}$ for $s,t \in I$.

Next we can find $\Phi_n$ (for $n < \omega$) such that:
\mn
\begin{enumerate}
\item[$\boxplus$]   $(a) \quad \Phi_n$ is 
proper for linear order and $\Phi_0 = \Phi$
\sn
\item[${{}}$]   $(b) \quad$ EM$_{\tau(\Phi)}(I,\Phi_n) \prec 
\text{ EM}_{\tau(\Phi)}(I,\Phi_{n+1})$ for every linear order $I$ and
$n < \omega$; 

\hskip25pt moreover
\sn
\item[${{}}$]   $(b)^+ \quad \tau(\Phi_n) \subseteq \tau(\Phi_{n+1})$ and
EM$(I,\Phi_n) \prec \text{ EM}_{\tau(\Phi_n)}(I,\Phi_{n+1})$ for every

\hskip25pt $n < \omega$ and linear order $I$
\sn
\item[${{}}$]   $(c) \quad$ if $|I| \le n$ then EM$_{\tau(\Phi)}(I,\Phi_n) = 
\text{ EM}_{\tau(\Phi)}(I,\Phi_{n+1})$ and 

\hskip25pt EM$_{\tau(T)}(I,\Phi_n) \cong M^*$
\sn
\item[${{}}$]   $(d) \quad |\tau(\Phi_n)| = \lambda$.
\end{enumerate}
\mn
This is easy.  Let $\Phi_\omega$ be the limit of $\langle \Phi_n:n < \omega
\rangle$, i.e. $\tau(\Phi_\omega) = \cup\{\tau(\Phi_n):n < \omega\}$
and if $k < \omega$ then $\EM_{\tau(\Phi_k)}(I,\Phi_\omega) =
\cup\{\EM_{\tau(\Phi_k)}(I,\Phi_n):n \in [k,\omega)\}$.  
So as $M^*$ is a superlimit model, for any linear order $I$
of cardinality $\lambda,\EM_{\tau(T)}(I,\Phi_\omega)$ is  the direct
limit of $\langle \EM_{\tau(T)}(J,\Phi_\omega):J \subseteq I$ 
finite$\rangle$, each isomorphic to $M^*$, so as we have assumed that
$M^*$ is a superlimit model it follows that EM$_{\tau(T)}(I,\Phi_\omega)$
is isomorphic to $M^*$.  But by \cite[III]{Sh:300} or \cite{Sh:E59}
which may eventually be \cite[III]{Sh:e} there are $2^\lambda$ many pairwise
non-isomorphic models of this form varying $I$ on the linear orders of
cardinality $\lambda$, contradiction.

\noindent
2) First assume $\lambda = \lambda^{\aleph_0}$.  Let $\tau \subseteq
\tau_T$ be countable such that $T' = T \cap \bbL(\tau)$ is not
superstable.  Clearly if $M^*$ is $(\lambda,\aleph_0)$-limit model then $M^*
\restriction \tau'$ is not $\aleph_1$-saturated.

\noindent
[Why?  As in \cite[Ch.VI,\S6]{Sh:a}, but we shall give full details.  
There are $N_* \models T,p
  = \{\varphi_n(\lambda,\bar a_n):n < \omega\}$ a type in $N_*,\bar
  a_n \triangleleft \bar a_{n+1},\bar a_{<>}$ empty and
  $\varphi_{n+1}(x,\bar a_{n+1})$ forks over $\bar a_n$.  Let $\bold
  F(M)$ be such that if $n < \omega$ and $\bar b_n \subseteq M$
  realizes $\tp(\bar a_n,\emptyset,N_*)$ then for some $\bar b_{n+1}$
  from $\bold F,M$ realizing $\tp(\bar a_{n+1},\emptyset,N_*)$, the
  type $\tp(\bar b_{n+1},M,\bold F(M))$ does not fork over $b_n$.]
But if $\kappa =
\text{ cf}(\kappa) \in [\aleph_1,\lambda]$ and $M^*$ is
a $(\lambda,\kappa)$-limit then $M^* \restriction \tau'$ is
$\aleph_1$-saturated, contradiction.

The case $\lambda \ge |T| + \beth_\omega$ is more complicated (the
assumption $\lambda \ge \beth_\omega$ is to enable us to use
\cite{Sh:460} or see \cite{Sh:829} for a simpler proof; we can 
use weaker but less transparent assumptions;
maybe $\lambda \ge 2^{\aleph_0}$ suffices).  

As $T$ is stable not superstable by \cite{Sh:c} for some $\bar\Delta$:
\mn
\begin{enumerate}
\item[$\circledast_1$]    for any $\mu$ there are $M$ and 
$\langle a_{\eta,\alpha}:\eta \in {}^\omega \mu$ and
$\alpha < \mu \rangle$ such that
\begin{enumerate}
\item[$(a)$]   $M$ is a model of $T$  
\sn
\item[$(b)$]   $\bold I_\eta = \{a_{\eta,\alpha}:\alpha < \mu\}
\subseteq M$ is an indiscernible set (and $\alpha < \beta < \mu
\Rightarrow a_{\eta,\alpha} \ne a_{\eta,\beta}$)
\sn
\item[$(c)$]   $\bar\Delta = \langle \Delta_n:n < \omega
\rangle$ and $\Delta_n \subseteq \Bbb L_{\tau(T)}$ infinite
\sn
\item[$(d)$]   for $\eta,\nu \in {}^\omega \mu$ we have
Av$_{\Delta_n}(M,\bold I_\eta) = \text{ Av}_{\Delta_n}(M,\bold I_\nu)$
iff $\eta \restriction n = \nu \restriction n$.
\end{enumerate}
\end{enumerate}
\mn
Hence by \cite[VIII]{Sh:c}, or see \cite{Sh:E59} 
assuming $M^*$ is a universal model of $T$
of cardinality $\lambda$ \,:
\mn
\begin{enumerate}
\item[$\circledast_{2.1}$]   there is $\Phi$ such that
\begin{enumerate}
\item[$(a)$]   $\Phi$ is proper for $K^\omega_{\text{tr}},\tau_T
\subseteq \tau(\Phi),|\tau(\Phi)| = \lambda \ge |T| + \aleph_0$
\sn
\item[$(b)$]   for $I \subseteq {}^{\omega \ge} \lambda$,
EM$_{\tau(\Phi)}(I,\Phi)$ is a model of $T$ and $I \subseteq J
\Rightarrow \text{ EM}(I,\Phi) \prec \text{ EM}(J,\Phi)$
\sn
\item[$(c)$]   for some two-place function symbol $F$ if for $I \in
K^\omega_{\text{tr}}$ and $\eta \in P^I_\omega,I$ a subtree of
${}^{\omega \ge}\lambda$ for transparency we let 
$\bold I_{I,\eta} = \{F(a_\eta,a_\nu):\nu \in I\}$ then $\langle \bold
I_{I,\eta}:\eta \in P^I_\omega\rangle$ are as in
$\circledast_1(b),(d)$.
\end{enumerate}
\end{enumerate}
\mn
Also
\mn
\begin{enumerate}
\item[$\circledast_{2.2}$]   if $\Phi_1$ satisfies (a),(b),(c) of
$\circledast_{2.1}$ and $M$ is a universal model of $T$ then there is
$\Phi^*_2$ satisfying (a),(b),(c) of $\circledast_{2.1}$ and $\Phi_1 \le
\Phi^*_2$ see $\circledast_{2.3}(a)$ and for every finitely generated $J
\in K^\omega_{\text{tr}}$, see $\circledast_{2.3}(b)$
below, there is $M' \cong M$ such that
$\EM_{\tau(T)},(J,\Phi_1) \prec M' \prec \EM_{\tau(T)}(J,\Phi^*_2)$
\sn
\item[$\circledast_{2.3}$]   $(a) \quad$ we say
$\Phi_1 \le \Phi_2$ when $\tau(\Phi_1) \subseteq \tau(\Phi_2)$ and
$J \in K^\omega_{\tr} \Rightarrow$

\hskip25pt $\EM(J,\Phi_1) \prec \EM_{\tau(\Phi_1)}(J,\Phi_2)$
\sn
\item[${{}}$]   $(b) \quad$ we say $J \subseteq I$ is finitely generated if it
has the form $\{\eta_\ell:\ell < n\} \cup$ 

\hskip25pt $\{\rho$: for some $n,\ell$ we have 
$\rho \in P^I_n$ and $\rho <^I \eta_\ell\}$
for some 

\hskip25pt  $\eta_0,\dotsc,\eta_{n-1} \in P^I_\omega$
\sn
\item[$\circledast_{2.4}$]   if $M_* \in \text{ EC}_\lambda(T)$ is
superlimit (or just weakly $S$-limit, $S \subseteq \lambda^+$
stationary) \then \, there is $\Phi$ as in $\circledast_{2.1}$ above
such that EM$_{\tau(T)}(J,\Phi) \cong M_*$ for every finitely
generated $J \in K^\omega_{\text{tr}}$
\sn
\item[$\circledast_{2.5}$]   we fix $\Phi$ as in
$\circledast_{2.4}$ for $M_* \in \text{ EC}_\lambda(T)$ superlimit.
\end{enumerate}
\mn
Hence (mainly by clause (b) of $\circledast_{2.1}$ and
$\circledast_{2.4}$ as in the proof of part (1))
\mn
\begin{enumerate}
\item[$\circledast_3$]   if $I \in K^\omega_{\text{tr}}$ has
cardinality $\le \lambda$ then EM$_{\tau(\Phi)}(I,\Phi)$ is isomorphic
to $M^*$.
\end{enumerate}
\mn
Now by \cite{Sh:460}, we can find regular uncountable $\kappa <
\beth_\omega$ such that $\lambda = \lambda^{[\kappa]}$, see Definition
\ref{c8}. 

Let $S = \{\delta < \kappa:\text{cf}(\delta) = \aleph_0\}$ and $\bar \eta =
\langle \eta_\delta:\delta \in S\rangle$ be such that 
$\eta_\delta$ an increasing
sequence of length $\omega$ with limit $\delta$.

For a model $M$ of $T$ let OB$_{\bar \eta}(M) = \{\bar{\bold
a}:\bar{\bold a} = \langle a_{\eta_\delta,\alpha}:\delta \in W$ 
and $\alpha < \kappa\rangle,W
\subseteq S$ and in $M$ they are as in $\circledast_1(b),(d)\}$.

\noindent
For $\bar{\bold a} \in \text{ OB}_{\bar \eta}(M)$ let $W[\bar{\bold
a}]$ be $W$ as above and let

\begin{equation*}
\begin{array}{clcr}
\Xi(\bar{\bold a},M) = \{\eta \in {}^\omega \kappa:&\text{ there is an
indiscernible set} \\
  &\bold I = \{a_\alpha:\alpha < \kappa\} \text{ in } M \text{ such
  that for every } n \\
  &\text{ for some } \delta \in W[\bar{\bold a}],\eta \restriction n =
  \eta_\delta \restriction n \text{ and} \\
  &\text{ Av}_{\Delta_n}(M,\bold I) = \text{
  Av}_{\Delta_n}(M,\{a_{\eta_\delta,\alpha}:\alpha < \kappa\})\}.
\end{array}
\end{equation*}
\mn
Clearly
\mn
\begin{enumerate}
\item[$\circledast_4$]   $(a) \quad$ if $M \prec N$ then OB$_{\bar
\eta}(M)  \subseteq \text{ OB}_{\bar \eta}(N)$
\sn
\item[${{}}$]   $(b) \quad$ if $M \prec N$ and $\bar{\bold a} \in
\text{ OB}_{\bar \eta}(M)$ then $\Xi(\bar{\bold a},M) \subseteq
\Xi(\bar{\bold a},N)$.
\end{enumerate}
\mn
Now by the choice of $\kappa$ it should be clear that
\mn
\begin{enumerate} 
\item[$\circledast_5$]   if $M \models T$ is of cardinality
$\lambda$ then we can find an elementary extension $N$ of $M$ of
cardinality $\lambda$ such that for every $\bar{\bold a} \in \text{
OB}_{\bar\eta}(M)$ with $W[\bar{\bold a}]$ a stationary subset of
$\kappa$, for some stationary $W' \subseteq W[\bar{\bold a}]$ the set 
$\Xi[\bar{\bold a},N]$ includes $\{\eta \in {}^\omega \kappa:
(\forall n)(\exists
\delta \in W')(\eta \restriction n = \eta_\delta \restriction n)\}$,
(moreover we can even find $\varepsilon^* < \kappa$ and 
$W_\varepsilon \subseteq W$ for $\varepsilon < \varepsilon^*$ satisfying
$W[\bar{\bold a}] = \cup\{W_\varepsilon:\varepsilon < \varepsilon^*\}$)
\sn
\item[$\circledast_6$]   we can find $M \in \EC_\lambda(T)$
isomorphic to $M^*$ such that for every $\bar{\bold a} \in 
\OB_{\bar\eta}(M)$ with $W[\bar{\bold a}]$ a stationary subset of
$\kappa$, we can find a stationary subset $W'$ of $W[\bar{\bold a}]$
such that the set $\Xi[\bar{\bold a},M]$ includes
$\{\eta \in {}^\omega \mu:(\forall n)(\exists \delta \in W')
(\eta \restriction n = \eta_\delta \restriction n)\}$.
\end{enumerate}
\mn
[Why?  We choose $(M_i,N_i)$ for $i < \kappa^+$ such that
\mn
\begin{enumerate}
\item[$(a)$]   $M_i \in \text{ EC}_\lambda(T)$ is
$\prec$-increasing continuous
\sn
\item[$(a)$]   $M_{i+1}$ is isomorphic to $M^*$
\sn
\item[$(a)$]  $M_i \prec N_i \prec M_{i+1}$
\sn
\item[$(a)$]  $(M_i,N_i)$ are like $(M,N)$ in $\circledast_5$.
\end{enumerate}
\mn
Now $M = \cup\{M_i:i < \kappa^+\}$ is as required.

\noindent
Now the model $M$ is isomorphic to $M^*$ as $M^*$ is superlimit.]

Now the model from $\circledast_6$ is not isomorphic to $M' = \text{
EM}_{\tau(T)}({}^{\omega >}\lambda \cup \{\eta_\delta:\delta \in
S\},\Phi)$ where $\Phi$ is from $\circledast_{2.1}$.  But $M' \cong M^*$
by $\circledast_3$.

\noindent
Together we are done.  
\end{PROOF}

\noindent
The following claim says in particular that
if some not unreasonable {\rm pcf} conjectures holds, the conclusion 
holds for every $\lambda \ge 2^{\aleph_0}$.
\begin{claim}
\label{nlm.0.1g}  
Assume $T$ is stable not superstable, $\lambda \ge |T|$ and 
$\lambda \ge \kappa = \text{\rm cf}(\kappa) > \aleph_0$.

\noindent
1) $T$ has no $(\lambda,\kappa)$-superlimit model provided that
$\kappa = \cf(\kappa) > \aleph_0,\lambda \ge 
\kappa^{\aleph_0}$ and $\lambda = \bold U_D(\lambda) := 
\Min\{|{\cP}|:{\cP} \subseteq [\lambda]^\kappa$ and for every
$f:\kappa \rightarrow \lambda$ for some $u \in {\cP}$ we have
$\{\alpha < \kappa:f(\alpha) \in u\} \in D^+$, where $D$ is a normal
filter on $\kappa$ to which $\{\delta < \kappa:\cf(\delta) =
\aleph_0\}$ belongs.

\noindent
2) Similarly if $\lambda \ge 2^{\aleph_0}$ and letting $J_0 = \{u
   \subseteq \kappa:|u| \le \aleph_0\},J_1 = \{u \subseteq \kappa:u
   \cap S^\kappa_{\aleph_0}$ non-stationary$\}$ we have $\lambda =
   \bold U_{J_1,J_0}(\lambda) := \text{\rm Min}\{|\cP|:\cP \subseteq
   [\lambda]^{\aleph_0}$, if $u \in J_1,f:(\kappa \backslash u)
   \rightarrow \lambda$ then for some countable infinite $w \subseteq
   \kappa(u)$ and $v \in \cP$, {\rm Rang}$(f \rest w) \subseteq v\}$.
\end{claim}

\begin{proof}
Like \ref{nlm.0.1}(2).
\end{proof}

\begin{claim}
\label{n9}
1) Assume $T$ is unstable and $\lambda \ge |T| + \beth_\omega$.  \Then \,
   for at most one regular $\kappa \le \lambda$ does $T$ have a weakly
   $(\lambda,\kappa)$-limit model and even a weakly
   $(\lambda,S)$-limit model for some stationary $S \subseteq
   S^\lambda_\kappa$.

\noindent
2) Assume $T$ is unsuperstable and $\lambda \ge |T| +
   \beth_\omega(\kappa_2)$ and $\kappa_1 = \aleph_0 < \kappa_2 =
   \cf(\kappa_2)$.  Then $T$ has no model which is a weak
   $(\lambda,S)$-limit where $S \subseteq \lambda$ and $S \cap
   S^\lambda_{\kappa_\ell}$ is stationary for $\ell=1,2$.
\end{claim}

\begin{PROOF}{\ref{n9}}
1) Assume $\kappa_1 \ne \kappa_2$ form a counterexample.  Let $\kappa
   < \beth_\omega$ be regular large enough such that $\lambda =
   \lambda^{[\kappa]}$, see Definition \ref{c8} and $\kappa
   \notin \{\kappa_1,\kappa_2\}$.  Let $m,\varphi(\bar x,\bar y)$ be
   as in the proof of \ref{nlm.0.1}
\mn
\begin{enumerate}
\item[$(*)$]  if $M \in \EC_\lambda(T)$ then there is $N$ such
that
\begin{enumerate}
\item[$(a)$]  $N \in \EC_\lambda(T)$
\sn
\item[$(b)$]  $M \prec N$
\sn
\item[$(c)$]  if $\bar{\bold a} = \langle \bar a_i:i < \kappa \rangle
\in {}^\kappa({}^m M)$ for $\alpha < \kappa$ then for some $\cU \in
[\kappa]^\chi$ for every uniform ultrafilter $D$ on $\kappa$ to which
$\cU$ belongs there is $\bar a_D \in {}^n N$ such that $\tp(\bar a_D,
N,N) = \Av(\bar{\bold a}/D,M) = \{\psi(\bar x,\bar c):
\psi(\bar x,\bar z) \in \bbL(\tau_T),\bar c \in {}^{\ell g(\bar
z)}M$ and $\{\{\alpha < \kappa:N \models \psi[\bar a_{i_\alpha},\bar c]\} 
\in D\}$.
\end{enumerate}
\end{enumerate}
\mn
Similarly
\mn
\begin{enumerate}
\item[$\boxplus_1$]  for every function $\bold F$ with domain $\{\bar
M:\bar M$ an $\prec$-increasing sequence of models of $T$ of length $<
\lambda^+$ each with universe $\in \lambda^+\}$ such that $M_i \prec
\bold F(\bar M)$ for $i < \ell g(\bar M)$ and $\bold F(\bar M)$ has
universe $\in \lambda^+$ \underline{there} is a sequence $\langle
M_\varepsilon:\varepsilon < \lambda^+\rangle$ obeying $\bold F$ such
that: for every $\varepsilon < \lambda^+$ and $\bar{\bold a} \in
{}^\kappa({}^m(M_\varepsilon))$ for $\alpha < \kappa$, there is $\cU
\in [\kappa]^\kappa$ such that for every ultrafilter $D$ on $\kappa$
to which $\cU$ belongs, for every $\zeta \in (\varepsilon,\lambda^+)$
there is $\bar{\bold a}_{D,\zeta} \in {}^m(M_{\zeta +1})$ realizing
$\Av(\bar{\bold a}/D,M_\zeta)$ in $M_{\zeta +1}$.
\end{enumerate}
\mn
Hence
\mn
\begin{enumerate}
\item[$\boxplus_2$]  for $\langle M_\alpha:\alpha < \lambda^+\rangle$
as in $\boxplus_1$ for every limit $\delta < \lambda^+$ of cofinality
$\ne \kappa$ for every $\bar{\bold a} = \langle \bar a_i:i <
\kappa\rangle \in {}^\kappa({}^m(M_\delta))$, there is $\cU \in
[\kappa]^\kappa$ such that for every ultrafilter $D$ on $\kappa$
to which $\cU$ belongs, there is a sequence $\langle \bar
b_\varepsilon:\varepsilon < \cf(\delta)\rangle \in
{}^{\text{cf}(\delta)}({}^m(M_\delta))$ such that for every $\psi(\bar
x,\bar z) \in \bbL(\tau_T)$ and $\bar c \in {}^{\ell g(\bar
z)}(M_\delta)$ for every $\varepsilon < \cf(\delta)$ large
enough, $M_\delta \models \psi[\bar b_\varepsilon,\bar c]$ iff
$\psi(\bar x,\bar c) \in \Av(\bar{\bold a}/D,M_\delta)$.
\end{enumerate}
\mn
The rest should be clear.

\noindent
2) Combine the above and the proof of \ref{nlm.0.1}(2).
\end{PROOF}
\newpage

\bibliographystyle{alphacolon}
\bibliography{lista,listb,listx,listf,liste,listz}

\end{document}